\documentclass[11pt]{article}
\usepackage{amssymb,amsmath}
\newtheorem{theorem}{Theorem}
\newtheorem{corollary}{Corollary}[section]
\newtheorem{lemma}[corollary]{Lemma}
\newtheorem{proposition}[corollary]{Proposition}

\newcommand{\Prob} {{\bf P}}
\newcommand{\Z}{{\mathbb Z}} 
\newcommand{\E}{{\bf E}}
\newcommand{\Es}{{\rm Es}}

\newcommand{\R}{{\mathbb R}}
\newcommand{\C}{{\mathbb C}}

\newcommand{\dist}{{\rm dist}}

\newcommand{\cA}{{\cal A}}
\newcommand{\disk}{{\cal B}}
\newcommand{\lf}{\lfloor}
\newcommand{\rf}{\rfloor}

\begin{document}

$\;$ \vspace{12ex}
\begin{center}
{\Large A Lower Bound on the Growth Exponent for \\
 Loop-Erased Random Walk in Two Dimensions}

\vspace{6ex}
Gregory F. Lawler \footnote{Research supported by the National Science Foundation}\\
Department of Mathematics\\
Box 90320\\
Duke University\\
Durham, NC 27708-0320\\

\vspace{2ex}

\end{center}

\begin{abstract}
The growth exponent $\alpha$ for loop-erased or Laplacian random walk
on the integer lattice is defined by saying that the expected time to
reach the sphere of radius $n$ is of order $n^\alpha$.  We prove that
in two dimensions, the growth exponent is strictly greater than one.
The proof uses a known estimate on the third moment of the escape
probability and an improvement on the discrete Beurling projection theorem.
\end{abstract}

\section{Introduction}

Loop-erased or Laplacian random walk (LERW) on the integer lattice $\Z^d, d \geq 2$ is a nonMarkovian, nearest neighbor, self-avoiding process.  This process was originally studied because it is a nontrivial process
that is self-avoiding, although there is very strong numerical and nonrigorous analytic evidence \cite{Duplantier,GB,book,Majumdar}
to believe that it is not in the same universality class as the usual self-avoiding walk.  There has been a recent interest in the LERW because of the connection between loop-erased walk and uniform spanning trees
\cite{Majumdar,robin,PW}.  

There are a number of ways to define the loop-erased walk.  If $d \geq 3$, one can take an infinite simple random walk and erase the loops chronologically to produce a self-avoiding path.  This is well-defined for $d \geq 3$ since the simple random walk is transient.  A somewhat modified definition is needed in two dimensions, and since this paper will concentrate on $d=2$ we will give this definition, which also works for higher dimensions.  If $\omega = [\omega(0),
\ldots,\omega(n)]$ is any nearest neighbor path in $\Z^d$, we define the loop-erased
path $L(\omega)$ as follows.  Let
\[  s_0 = \sup\{k \leq n: \omega(k) = \omega(0) \}, \]
and if $s_j < n$,
\[  s_{j+1} = \sup\{k \leq n:  \omega(k) = \omega(s_j + 1) \} . \]
If $l$ is the smallest index so that $s_l = n$, the loop-erased path
is given by
\[  L(\omega) = L[\omega(0),\ldots,\omega(n)] = [\omega(s_0),\ldots,\omega(s_l)].\]
Note the $L(\omega)$ is a self-avoiding,
 nearest neighbor path whose initial and end
points are the same as those of $\omega$. 

Let $S(k)$ denote a simple random walk in $\Z^d, d \geq 2$, starting at the
origin.  Let
\[  	C_m = \{x \in \Z^d: |x| < m \} , \]
with boundary
\[  \partial C_m = \{x \in \Z^d \setminus C_m: |y-x| = 1 \mbox{ for some }
      y \in C_m \} . \]
Let
\[  \sigma_m = \inf\{k: S(k) \in \partial C_m \} . \]
Let $\Lambda_n$ be the set of nearest neighbor, self-avoiding paths
$\omega = [0=\omega(0),\ldots,\omega(l)]$ of any length $l$ such that
$\omega(k) \in C_n$, $k < l$, and $\omega(l) \in \partial C_n$.  If $n \leq m$,
there is a measure on $\Lambda_n$, $\mu_{n.m}$ obtained by  considering the unique initial segment of $L(S[0,\sigma_m])$
that is in $\Lambda_n$.  We define the measure $\mu_n$ as the limit as $m$ tends
to infinity of the measures $\mu_{n,m}$.  It can be shown \cite[Chapter 7]{book}
 that the limit exists,
and the measures $\{\mu_n\}$ are consistent.  Hence this gives a measure on infinite
self-avoiding paths.  This measure is the same as the measure produced by the following
nonMarkovian transition probabilities.  If $x \in \Z^d$, and $\omega
=[\omega(0),\ldots,\omega(k)]$ is a finite
self-avoiding path, let
\[  f_m(x,\omega) = \Prob^x\{S[0,\sigma_m] \cap \omega = \emptyset \} \]
(here we make the natural identification of a path with its range).
Then if $|x-\omega(k)| = 1$,
\[  \Prob\{\hat{S}(k+1) = x \mid [\hat{S}(0),\ldots,\hat{S}(k)] = \omega \} =  
    \lim_{m \rightarrow \infty} \frac{f_m(x,\omega)}{\sum_{|y-\omega(k)| = 1}
      f_m(y,\omega)} . \]
Also from the discrete Harnack principle \cite[Theorem 1.7.6]{book},
 there exists a constant $c$ such that
if $\omega \subset C^{m/2}$,
\[ c^{-1}   \frac{f_m(x,\omega)}{\sum_{|y-\omega(k)| = 1}
      f_m(y,\omega)}  \leq \Prob\{\hat{S}(k+1) = x \mid [\hat{S}(0),\ldots,\hat{S}(k)] = \omega \} \leq c \frac{f_m(x,\omega)}{\sum_{|y-\omega(k)| = 1}f_m(y,\omega)} . \]
Hence, if $n \leq m/2$, $\omega \in \Lambda_n$,
\begin{equation}  \label{nov24.1}
 c^{-1}  \mu_{n} \leq \mu_{n,m} \leq c \mu_n . 
\end{equation}

We are interested in the exponent that measures the rate of growth
of $\hat{S}$.  This exponent is often phrased in terms of the
mean squared distance $\E[|\hat{S}(n)|^2]$.  We will use a different, but presumably
equivalent, formulation.  Let
\[ \hat{\sigma}_n = \inf\{k: \hat{S}(k) \in \partial C_n \} . \]
We define the exponent $\alpha = \alpha_d$ by the relation
\[  \E[\hat{\sigma}_n] \approx n^{\alpha} . \]
Intuitively we say that the paths of the loop-erased walk have fractal dimension
$\alpha$.  If $d \geq 4$, it is known \cite{book,Kahane}
that $\alpha = 2$,
with logarithmic corrections in four dimensions.
We have  no proof that $\alpha$ exists for $d =2,3$, 
so to be precise we should 
 define $\underline{\alpha}$ and $\overline{\alpha}$ to be the
$\liminf$ and $\limsup$, respectively, of 
\[ \frac{\log \E[\hat{\sigma}_n]}{\log n} . \]
By slight abuse of notation we will write $\alpha \geq s$ to mean
$\underline{\alpha} \geq s $ and $\alpha \leq s$ for $\overline{\alpha}
\leq s$.   For $d=2,3$, there
is an upper bound \cite{book}
\begin{equation}  \label{0.0}
 \alpha_d \leq \frac{d + 2}{3} , 
\end{equation}
but numerical simulations \cite{GB} indicate
that this bound is not sharp.  The right hand
side is exactly the Flory predictions for the corresponding quantity
for the usual self-avoiding walk.  This prediction is still expected
to be correct for the usual self-avoiding walk in two dimensions, but
is expected to be slightly lower than the actual value in three
dimensions \cite{MS}.  The conjectures imply that the LERW goes
to infinity faster than the usual simple random walk. 
By comparison to uniform spanning trees and the Potts model, nonrigorous
conformal field theory \cite{Duplantier,Majumdar} has been used to conjecture that
\[   \alpha_2 = \frac{5}{4} . \]
There is no reason to believe that $\alpha_3$ is a nice rational
number; numerical simulations do suggest that $\alpha_3 < 5/3$. 
The rigorous inequality
$\alpha \geq 1$ is immediate.  By comparison
to the intersection exponent in three dimensions \cite{JPhysA}, it can
be seen that $\alpha_3 > 1$, but no such inequality has
been shown in two dimensions.  The purpose of this paper is to 
prove 
\[  \alpha_2 > 1. \]	

\begin{theorem} \label{theorem.main}
There exist positive constants $c,\epsilon$ such that if $d=2$,
for all $n$,
\[   \E[\hat{\sigma}_n] \geq c n^{1 + \epsilon} . \]
\end{theorem}

Throughout this paper we use $c,c_1,c_2$ for positive constants whose
value may change from line to line.
For the remainder of the paper we will assume that $d=2$. By (\ref{nov24.1}),
to prove this estimate for a given $n$ it suffices to consider the
simple random walk up to time $\sigma_{2n}$; erase loops from this
path; and give an appropriate lower bound on
 the number of points before time $\sigma_n$ that are
not erased.  Let $V(j,n)$ be the event that $j \leq \sigma_{n}$ and that
$S(j)$ is not erased in producing $L(S[0,\sigma_{2n}])$. More
precisely,
\[  V(j,n) =
  \{j \leq \sigma_n; L(S[0,j]) \cap S[j+1,\sigma_{2n}] = \emptyset \} . \]
Note that if $\hat{\sigma}_{n,2n}$ is the analogue of $\hat{\sigma}_n$ for
walks stopped upon reaching $\partial C_{2n}$,
\[  \hat{\sigma}_{n,2n} \geq \sum_{j=0}^\infty I[V(j,n)] , \]
where $I$ denotes the indicator function.
Hence in order to prove the theorem it suffices to prove that there exist
$c,\epsilon$ such that for $n^2 \leq j \leq 2 n^2$,
\begin{equation}  \label{0.1}
  \Prob[V(j,n)] \geq c n^{-1 + \epsilon} . 
\end{equation}
The distribution of $L(S[0,j])$ is the same whether we consider
$S(0)$ or $S(j)$ as the origin, i.e. ``reverse'' loop-erased walk has the
same distribution as LERW \cite[Lemma 7.2.1]{book}.  To prove (\ref{0.1}) it
suffices to prove
the following estimate.  Let $S^1,S^2$ be independent simple random walks starting
at the origin with corresponding stopping times $\sigma^1_n,\sigma^2_n$.  Then
for some $c,\epsilon$,
\begin{equation}  \label{0.2}
  \Prob\{L(S^1[0,\sigma_n^1]) \cap S^2(0,\sigma^2_n] = \emptyset \}
   \geq c n^{-1 + \epsilon} . 
\end{equation}
That is, the probability that a simple walk and a loop-erased walk get to distance
$n$ without intersecting is greater than $c n^{-1+\epsilon}$.  This is the
main estimate of this paper.

Assume for ease that $S^i$ is defined on the probability space $(\Omega_i,
\Prob_i)$ and let $\E_i$ denote expectations with respect to
$\Prob_i$.  Define the $\Omega_1$ random variable
\[  X_n = \Prob_2\{S^2(0,\sigma^2_n] \cap L(S^1[0,\sigma^1_n]) = \emptyset \}. 
\]
It can be shown that,
\begin{equation}  \label{0.3}
  \E[X_n^3] \asymp n^{-2} , 
\end{equation}
where $\asymp$ means that both sides are comparable, i.e.,
bounded by a constant times a multiple
of the other side.  The upper bound for $\E[X_n^3]$, at least up to logarithmic
corrections, 
 can be found in \cite{book}; this and the inequality
$\E[X_n^3] \geq [\E(X_n)]^3$ give the estimate (\ref{0.0}).  The lower bound is what is
needed in this paper; the proof   will appear in \cite{Chad} as well
as the analogous result for three dimensions, but since we need it, we
will give a proof of
 \begin{equation}  \label{dec30}
\E[X_n^3] \geq c n^{-2}
\end{equation}
in this paper.    

The discrete Beurling projection theorem \cite{Kesten}
states that if
$\omega$ is any simple random walk path
in $\Z^2$ connecting the origin with $\partial C_n$,
then
\[  \Prob_2\{S^2(0,\sigma^2_n] \cap \omega = \emptyset \} \leq c n^{-1/2} . \]
The probability on the left is maximized (at least up to a multiplicative constant)
when  $\omega$ is a half line, in which case
 this probability is also greater than $c_1 
n^{-1/2}$.  This theorem
 implies that $X_n \leq c n^{-1/2}$, and hence 
\[  \E[X_n] \geq c n \E[X_n^3] . \]
This is not good enough to get (\ref{0.2}); in fact this only allows us to conclude
the trivial inequality $\alpha \geq 1$.  However, a typical loop-erased walk is
more crooked than a straight line.  The main technical tool in this paper is
an upper bound on the escape probability for walks in terms of  a particular
quantity that measures crookedness.  In particular, we will  show that
there is a $c$ and a $\delta > 0$ such that for all large $n$
\begin{equation}  \label{0.4}
 \Prob_1\{X_n \geq  c n^{-1/2-\delta} \} \leq n^{-3}.
\end{equation}
In particular,
\[  \E_1[X_n^3; X_n \leq c n^{-1/2-\delta}] \geq c_1 n^{-2} ,  \]
and hence,
\[  \E_1[X_n] \geq c n^{-1 + 2 \delta} . \]
This gives (\ref{0.2}).

The improvement on the discrete
Beurling projection theorem is of independent interest.
Let $C^m = C_{e^m}$ and suppose $\omega$ is a self-avoiding
 random walk path starting
at the origin ending at $C^m$.  Let 
\[ \eta_k = \inf\{t: \omega(t) \in \partial C^k \}. \]
For each $0 < \delta < \pi/2$, let
\[    Z_{\delta,m} = Z_{\delta,m}(\omega) 
   = \sum_{k=0}^{m-1} I\{|\arg(\omega(\eta_{k+1})) - \arg
   (\omega(\eta_k))| \geq \delta \} . \]
We will prove the following.  If $S$ is a simple random walk, let
\[  \tau_m = \sigma^m = \inf\{t: S(t) \in \partial C^m \} . \]

\begin{proposition}  \label{prop.1}
For every $\delta,a >0$, there exist $c < \infty, \beta > 1/2$,
such that if $\omega$ is any self-avoiding
random walk path connecting $0$ with $\partial C^m$ with
\[    Z_{\delta,m} \geq a m , \]
then
\[  \Prob\{S(0,\tau_m] \cap \omega = \emptyset \}
   \leq c e^{-\beta m } . \]
\end{proposition}

We first prove a corresponding proposition for Brownian motion using
the relationship between harmonic measure and extremal distance.  The
result for random walk is obtained using a strong approximation derived
from the Skorohod embedding.

The outline of the paper is as follows.  In Section \ref{thirdsec}
 we derive
the estimate on $\E[X_n^3]$; the argument is similar to that in
\cite{Kahane} where a corresponding result was proved for $d=4$.
The next section gives exponential estimates on the probability
that the loop-erased walk is very straight.  The goal of
Section \ref{brownsec} is to 
prove the analogue of Proposition \ref{prop.1} for Brownian motion.  The 
derivation
uses a relationship between escape probabilities in two dimensions and
a quantity, extremal distance, of a domain.  The last two sections derive
the result for random walk.  Section \ref{strongsec} reviews the necessary
facts about the strong approximation and then Section \ref{walksec} uses
this to obtain the result.  A similar idea with added complications
 has been used in \cite{walkcut,
Emily} to estimate probabilities that paths of
random walks do not intersect each other
with the corresponding probabilities for Brownian motion.
 
\section{Third Moment}  \label{thirdsec}

In this section we derive (\ref{dec30}). If $A \subset \Z^2$, let
\[  \Es_n(A) = \Prob\{S(0,\sigma_n] \cap A = \emptyset \} . \]

\begin{lemma} \label{lemma.mom.1}
Suppose $A \subset \Z^2$ contains 
at least one nearest neighbor of the origin. Then for all $n \geq 1$,
\[  \Es_n(A \setminus \{0\}) \leq 4 \Es_n(A). \]
\end{lemma}

{\bf Proof.}
 Let $A^\prime = A \setminus \{0\}$, and 
\[   \eta_n = \inf\{t > 0: S(t) \in \partial C_n \cup \{0\} \; \} . \]
Since $A^\prime$ contains a nearest neighbor of the origin,
\[ \Prob\{S(0,\eta_n] \cap A^\prime \neq \emptyset \}
   \geq \frac{1}{4} . \]
Let $Z_n$ be the number of visits to $0$ by $S$ before leaving $C_n$,
\[  Z_n = \sum_{t=0}^{\sigma_n} I\{S(t) = 0 \} . \]
Then by iterating above
\[  \Prob\{Z_n = k + 1, S(0,\sigma_n] \cap A^\prime = \emptyset\}
   \leq (\frac{3}{4})^k \Prob\{S(0,\sigma_n] \cap A = \emptyset \} . \]
Summing over $k$ gives the lemma.

\vspace{4ex}

If $S$ is a simple random walk starting at the origin, define the
random variable
\[  X_n = \Es_n(L(S[0,\sigma_n])) . \]
It is the goal of this section to show that 
\[   \E[X_n^3] \geq c n^{-2} , \]
for some constant $c$.  By the lemma, it suffices to show that
\begin{equation}  \label{moment.1}
 \E[(X^\prime_n)^3] \geq c n^{-2} , 
\end{equation}
where
\[ X_n^\prime = \Es_n(L(S[0,\sigma_n]) \setminus \{0\}) . \]
Let $S^1,\ldots,S^4$ be independent simple random walks starting at
the origin with stopping times $\sigma_n^1,\ldots,\sigma_n^4$.
Let 
\[  \Theta_n = S^2(0,\sigma_n^2] \cup
    S^3(0,\sigma_n^3] \cup S^4(0,\sigma_n^4] . \]
Then, by independence,
\[ \E[(X^\prime_n)^3] = \Prob\{ \Theta_n \cap
     [L(S[0,\sigma_n^1]) \setminus \{0\}] = \emptyset  \} .\]

Let $[S^i(t)]_1,[S^i(t)]_2$ denote the first and second
components of the random walk $S^i$.  Let $\cA^1_n,\ldots,\cA^4_n$,
denote the rectangles,
\[  \cA^1_n = \{(x,y) \in \Z^2: -2n \leq x \leq \frac{n}{4} ;
      -\frac{n}{4} \leq y \leq \frac{n}{4} \}, \]
\[  \cA^2_n = \{(x,y) \in \Z^2:   -\frac{n}{4} \leq x \leq 2n  ;
      -\frac{n}{4} \leq y \leq \frac{n}{4}\} , \]
\[  \cA^3_n = \{(x,y) \in \Z^2:   -\frac{n}{4} \leq x \leq \frac{n}{4}   ;
      -2n \leq y \leq \frac{n}{4} \}, \]
\[  \cA^4_n = \{(x,y) \in \Z^2:   -\frac{n}{4} \leq x \leq \frac{n}{4}   ;
      - \frac{n}{4}\leq y \leq 2n \} . \]
For $n^2 \leq t_1,t_2,t_3,t_4 \leq 2n^2$, let
$U_n(t_1,\ldots,t_4)$ be the event that the following holds:
\[      S^i[0,t_i] \subset \cA^i_n , \;\;\; i=1,2,3,4 ; \]
\[      |S^i(t)| \geq n , \;\;\; n^2 \leq t \leq t_i, \; i=1,2,3,4 . \]
Let $V_n(t_1,\ldots,t_4)$ be the event 
\[ V_n =  \{ \Theta_n \cap
     [L(S[0,\sigma_n^1]) \setminus \{0\}] = \emptyset  \}. \]
We will show that
\begin{equation}  \label{moment.2}
 \sum_{n^2 \leq t_1,t_2 \leq 2n^2} \Prob[U_n(t_1,t_2,n^2,n^2)]
    \cap V_n(t_1,t_2,n^2,n^2)] \geq c n^{2} . 
\end{equation}
By the obvious monotonicity this implies
\[   \Prob[U_n(n^2,n^2,n^2,n^2) \cap V_n(n^2,n^2,n^2,n^2)]
    \geq c n^{-2} . \]
From this (\ref{moment.1}) follows easily.

Let $S^5$ be a random walk starting at the origin and let
$G_n$ be the event that the following four conditions hold:
\[   [S^5(t)]_2 \leq \frac{n}{8} , \;\; 0 \leq t \leq 4n^2 , \]
\[    \frac{5n}{4} \leq [S^5(t)]_1  \leq \frac{7n}{4}, \;\;
           \frac{5}{4} n^2 \leq t \leq \frac{7}{4}n^2 \; \]
\[    [S^5(t)]_1 \leq n , \; \; 0 \leq t \leq n^2 , \; \]
\[    [S^5(t)]_1 \geq 3n, \;\; 2n^2 \leq t \leq 4n^2 . \]
It is easy to see from the invariance principle that there is a constant
$c$ such that
\[  \Prob[G_n] \geq c . \]
Now start another random walk $S^6$, independent of $S^5,$
starting at a point $x \in \Z^2$,
chosen from the discrete
ball of radius $n/16$ around $(3n/2,-2n)$.  Note that
the number of points in this ball is comparable to $n^2$.  Let
 $H_n$ be the event that the following conditions hold:
\[    \frac{11}{8}n \leq [S^6(t)]_1  \leq \frac{13}{8}n ,\;\;
     0 \leq t \leq 4 n^2 , \]
\[ -\frac{n}{4} \leq [S^6(t)]_2 \leq \frac{n}{4}, \;\;
           \frac{5}{4} n^2 \leq t \leq \frac{7}{4}n^2 \; \]
\[  [S^6(t)]_2 \leq -\frac{3n}{2}, \;\; 0 \leq t \leq n^2 , \]
\[  [S^6(t)]_2 \geq  \frac{3n}{2}, \;\;  2n^2 \leq t \leq 4 n^2 . \]
Again, by the invariance principle, there is a constant $c$, independent
of the starting point, such that if $x$ is in the ball of
radius $n/16$ around $(3n/2,-2n)$,
\[  \Prob^x[H_n] \geq c , \]
and hence for such $x$
 by independence $\Prob^x(G_n \cap H_n) \geq c $ (here the
$x$ refers to the starting point of the walk $H_n$). By summing, we get
\[  \sum_x \Prob^x(G_n \cap H_n) \geq c n^2. \]

Consider the loop-erased path $L(S^5[0,4n^2])$. On the event
$G_n \cap H_n$ we can see from geometric considerations that the paths
$L(S^5[0,4n^2])$ and $S^6[0,4n^2]$ must intersect, and the points
of intersection must occur in the set
\[  \{(x,y): \frac{11}{8}n \leq x \leq \frac{13}{8}n ,
           -\frac{1}{8}n \leq y \leq \frac{1}{8}n \} . \]
Let $\rho^1 = \rho^1_n$ be the smallest time t such that
the point $S^5(t)$ is included in the loop-erased path and $S^6[0,4n^2]$.
 More precisely, $\rho^1$ is defined by
the conditions:
\[  L(S^5[0,\rho^1]) \cap S^5[\rho^1 + 1,4n^2] 
   = \emptyset , \; \]
\[    S^5(\rho^1) \in S^6[0,4n^2] . \]
Note that on $G_n \cap H_n$, $n^2 \leq \rho^1 \leq 2n^2$.  Let
\[ \rho^2 = \inf\{t:  S^6(t) = S^5(\rho^1) \} . \]
Again, on $G_n \cap H_n, n^2 \leq \rho^2 \leq 2n^2$. 
Hence
\begin{equation}  \label{moment.3}
  \sum_{x } \sum_{n^2 \leq t_1,t_2,
    \leq 2n^2} \Prob^x[G_n \cap H_n; \rho_1 = t_1; \rho_2 = t_2 ]
    \geq c n^2 ,
\end{equation}
where the outer sum is over all $x$ in the ball of radius $n/16$ about
$(3n/2,-2n)$.

We define the random walks $S^1,\ldots,S^4$ in terms of $S^5,S^6$.
Let
\[  S^1(t) = S^5(\rho^1 - t) - S^5(\rho^1), \;\; t= 0,1,\ldots,\rho^1, \]
\[  S^2(t) = S^6(\rho^2 - t) - S^6(\rho^2), \;\; t=0,1,\ldots,\rho^2, \]
\[ S^3(t) = S^5(\rho^1 + t) - S^5(\rho^1), \;\; t=0,1,2,\ldots,n^2 , \]
\[  S^4(t) = S^6(\rho^2 + t) - S^6(\rho^2), \;\; t=0,1,2,\ldots,n^2 . \]
Then in terms of these random walks, we can see that
\[  \sum_{x} \Prob^x[G_n \cap H_n; \rho_1 = t_1; \rho_2 = t_2 ]
  \geq P[V(t_1,t_2,n^2,n^2)]. \]
This combined with (\ref{moment.3}) gives (\ref{moment.2}) and hence 
proves (\ref{moment.1}).

We note that a similar argument will work in $d=3$.  One complication
arises in that one must show that the loop-erasure of $S^5$ and
the random walk $S^6$ intersect with positive probability independent
of $n$.  This is not an issue in two dimensions since the continuous
curves must cross.  For a proof of
\[  \E[X_n^3] \geq c n^{-1} \]
in three dimensions, see \cite{Chad}.   This proof can also be adapted easily
to prove the upper bound
\[  \E[X_n^3] \leq c n^{d-4} . \]
Here, one essentially needs only the trivial estimate that the probability that
$S^5$ and $S^6$ intersect is of order one, although a few other minor
technicalities arise.

\section{Crookedness of Loop-Erased Walk}

Let $S$ be a simple random walk in $\Z^2$ and
\[ \tau_n = \sigma_{e^n} = \inf\{t: S(t) \in \partial C^n \} . \]
Let $\hat{S}_n(j)$ denote the walk obtained by erasing loops from
$S[0,\tau_n]$.  This gives a measure on self-avoiding paths starting
at the origin and ending upon reaching $\partial C^n$.  If $k \leq n$,
let
\[  \hat{\tau}_{k,n} = \inf\{t: \hat{S}_n(t) \in \partial C^k\} , \]
\[ \hat{\tau}_k = \inf\{t:\hat{S}(t) \in \partial C^k\} . \]
Let
\[  \hat{W}_{n,\delta} = \sum_{k=1}^n I\{|  \arg[\hat{S}_n(\hat{\tau}_{k,n})]
     - \arg   [\hat{S}_n(\hat{\tau}_{k-1,n})] \; | \leq \delta \} , \]
\[  \tilde{W}_{n,\delta} = \sum_{k=1}^n I\{|  \arg[\hat{S}(\hat{\tau}_{k})]
     - \arg   [\hat{S}(\hat{\tau}_{k-1})] \; | \leq \delta \} , \]

\begin{proposition}  \label{prop.crook.1}
For every $M < \infty$ and every $\epsilon > 0$, there exist $\delta >0$
and $c < \infty$ such that for all $n$,
\[   \Prob\{\hat{W}_{n,\delta} \geq \epsilon n\} \leq c e^{-Mn} . \]
\[ \Prob\{\tilde{W}_{n,\delta} \geq \epsilon n\} \leq c e^{-Mn} . \]
\end{proposition}

The second inequality follows immediately from the first and (\ref{nov24.1}),
so we will only prove the first.
We will actually prove the following stronger proposition.  Standard
large deviation estimates for binomial random variables can be used to
derive Proposition \ref{prop.crook.1} from Proposition \ref{prop.crook.2}.

\begin{proposition}  \label{prop.crook.2}
For every $\epsilon > 0$, there exist a $\delta > 0$ and a $K < \infty$,
such that
for any $K \leq k \leq n-1$,
\[  \Prob\{ |\arg[\hat{S}_n(\hat{\tau}_{k,n})]
     - \arg   [\hat{S}_n(\hat{\tau}_{k-1,n})] \; | \leq \delta \mid
   \hat{S}_n(t), t=0,\ldots,\hat{\tau}_{k-1,n}\}
   \leq \epsilon . \]
\end{proposition}

Fix $k \leq n-1$ and suppose we know
\[  [\hat{S}_n(0),\hat{S}_n(1),\ldots,\hat{S}_n(\hat{\tau}_{k-1,n})] . \]
Analysis of the loop-erasing procedure shows that to determine the
distribution
\[  \{\hat{S}_n(t): \hat{\tau}_{k-1,n} < t \leq \hat{\tau}_{n,n} \}, \]
we
can do the following:  start another simple random walk $S^1(t)$ at
$\hat{S}_n(\hat{\tau}_{k-1,n})$ with corresponding stopping times $\tau_n^1$
and condition the walk so that
\[  S^1(0,\tau^1_n] \cap \hat{S}[0,\hat{\tau}_{k-1,n}] = \emptyset ;\]
then, erase loops from this conditioned path.  Let $x \in \partial C^k$
and let ${\cal B}(x,\delta e^k)$ represent the discrete ball of radius
$\delta e^k$ about $x$.  Proposition \ref{prop.crook.3} will show that the
probability that the conditioned random walk enters this ball at any
time tends to zero as $\delta \rightarrow 0$.  Since the loop-erased
path is a subpath of the conditioned path, the probability it enters
the balls also tends to zero.  From this we conclude Proposition
\ref{prop.crook.2}. Before stating and proving this proposition, we
review some facts about simple random walks (see \cite[Section 1.6]{book}
for more details).  Suppose the simple random walk starts at the origin
and ${\cal B}$ is a discrete ball of radius $\delta e^n$ centered
at some $y \in \partial C_n$.  Let $ K_\delta$ be sufficiently large
such that
every discrete ball of radius $\delta e^{n}$ has
 has at least $\delta^2 e^{2n}$ points provided $n \geq K_\delta$.  If we let
$V$ denote the number of visits to ${\cal B}$ before leaving the ball
of radius $e^{n+1}$, then standard estimates give
\[  \E^0[V] \leq c_1 \delta^2 , \]
and if $z \in {\cal B}$, and $n \geq K_\delta$,
\[  \E^z[V] \geq c_2 \delta^2 \log(1/\delta). \]
In particular,
\[ \E^0[V \mid V \geq 1] \geq  c_2 \delta^2 \log(1/\delta). \]
and hence
\begin{equation}  \label{jan3}
\Prob^0\{ V \geq 1\} \leq c [\log(1/\delta)]^{-1}.
\end{equation}

\begin{proposition}  \label{prop.crook.3}
There exists a constant $c< \infty$  such that the following is true.
Let $K_\delta$ be  as above, $k < n-1$, and let $\omega = [\omega(0),\ldots,\omega(r)]$ be
a random walk path with $\omega(0) = 0; \omega(t) \in C^k, t < r;
\omega(k) \in \partial C^k$.  Let $S$ be a simple random walk with stopping
times $\tau_m$ and let ${\cal B}$ be a discrete ball of radius
$\delta e^k$ centered at  $y \in \partial C^{k+1}$.  Then
if $r \geq K_\delta$,
\[ \Prob^{\omega(k)}\{S(0,\tau_n] \cap {\cal B} \neq \emptyset 
\mid S(0,\tau_n] \cap \omega = \emptyset \} \leq c [\log (1/\delta) ]^{-1}. \]
\end{proposition}

{\bf Proof.} Fix $k_\delta \leq
k<n-1, \delta >0$, and $\omega,{\cal B}$ as in the statement
of the theorem (but constants in this proof are independent of all of
these).  Let
\[  a_1 = \inf\{t: S(t) \in \partial C^{k + (1/2)} \} , \]
\[  b_1 = \inf\{t: S(t) \in \partial C^{k+2} \}, \]
and for $j > 1$,
\[  a_j = \inf\{t > b_{j-1}: S(t) \in \partial C^{k + (1/2)} \} , \]
\[  d_j = \inf\{t > b_{j-1}: S(t) \in \partial C^{k + 1} \} , \]
\[  b_j = \inf\{t > a_{j-1}: S(t) \in \partial C^{k+2} \}. \]
Standard estimates (see \cite[Section 1.6]{book}and (\ref{jan3}))
tell us that for $k \geq K_\delta$,
$x \in \partial C^{k+2}$,
\[ \Prob^x\{ \tau_{k+1} > \tau_n \} \asymp
   \Prob^x\{ \tau_{k+(1/2)} > \tau_n \} \asymp c [n-k]^{-1} , \]
\[  \Prob\{S[a_j,b_j] \cap {\cal B} \neq \emptyset \mid
   S(t), t=0,1,\ldots,a_j\} \leq c [\log(1/\delta)]^{-1} , \]
\[  \Prob\{S[b_j,a_{j+1}] \cap {\cal B} \neq \emptyset \mid
   S(t), t=0,1,\ldots,b_j\} \leq c [\log(1/\delta)]^{-1} . \]
Also (see \cite[Lemma 2.5.3]{book}),
\[  \Prob\{S[a_j,b_j] \cap \omega \neq \emptyset \mid 
     S(t), t=0,1,\ldots,a_j\} \geq c . \]
Let $j^*$ be the largest $j$ such that
\[    \{b_j < \tau_n \} . \]
Then the estimates above show there exists a $u < 1$ such that
\[   \Prob\{j^* = j; S[0,b_j] \cap {\cal B} \neq \emptyset;
   S[0,b_j] \cap \omega = \emptyset \}
     \leq c j u^j [\log(1/\delta)]^{-1} [n-k]^{-1} . \]
Summing over all $j$ we get for $x \in \partial C^{k + (1/2)}$.
\[  \Prob^x\{S(0,\tau_n] \cap {\cal B} \neq \emptyset;
   S[0,b_j] \cap \omega = \emptyset \}
       \leq c [\log(1/\delta)]^{-1} [n-k]^{-1} . \]
But it is easy to see for $x \in \partial C^{k+(1/2)}$,
\[ \Prob^x\{S(0,\tau_n] \cap {\cal B} = \emptyset ; 
S[0,b_j] \cap \omega = \emptyset\}
   \geq c  [\log(1/\delta)]^{-1} [n-k]^{-1} . \]
(To see this we bound by the probability that 
 \[ S(0,a_1] \cap {\cal B} = \emptyset, S[0,b_j] \cap \omega = \emptyset \]
times the probability that $d_2 > \tau_n$.)  This gives the
proposition.

\section{Extremal Length and Escape Probabilities}  \label{brownsec}

There is a close relationship between escape probabilities for
planar Brownian motions and a quantity known as extremal length.
In this section we will review some of the basic facts about
extremal length, and then use an extremal length estimation to
estimate escape probabilities.  For more details about extremal
length, see \cite{Ahlfors,Pom}.

 Let $\Gamma$ be a collection
of piecewise smooth curves $\gamma: [0,1] \rightarrow \C$.
Suppose $\rho: \C
\rightarrow  \R$ is a measurable function with
\[    a(\rho) = \int_\C \rho^2 \; dx \; dy < \infty. \]
Let
\[ L(\rho) = L(\rho,\Gamma) = 
   \inf \int_\gamma \rho \; d|z| , \]
where the infimum is over all $\gamma \in \Gamma$.
 The  extremal length
$\Delta = \Delta(\Gamma)$
 is defined by
\begin{equation}  \label{extremal.1}
  \Delta = \sup \frac{L(\rho)^2}{a(\rho)} , 
\end{equation}
where the supremum is over all $\rho$ with $a(\rho) < \infty$.  
The case of greatest interest for us is where $D$ is a bounded
domain; $V_1,V_2$ closed subsets of $\C$; and $\Gamma = \Gamma(D,
V_1,V_2)$ is the set
of all piecewise smooth $\gamma$ with $\gamma(0) \in V_1,\gamma(1) \in
V_2$, and $\gamma(0,1) \subset D$.  In this case we write
$ \Delta = \Delta(D,V_1,V_2) $.
 Extremal length is a conformal invariant, i.e., if
$f:D \rightarrow D_1$ is a conformal transformation defined
up to the boundary then
\[   \Delta(f(D),f(V_1 \cap \bar{D}),f(V_2 \cap \bar{D}))
    = \Delta(D,V_1,V_2). \]

If $D$ is the rectangle
\[   D=\{x + iy: 0 < x < a, 0 < y < b \} ,\]
and
\[  V_1 = \{iy: 0 \leq y \leq b\}, \;\; V_2 = \{a +iy:
    0 \leq y \leq b \} , \]
then
\[    \Delta(D,V_1,V_2) = \frac{a}{b} , \]
with the supremum in (\ref{extremal.1}) being taken on by $\rho \equiv 1$.
If $D$ is the split annulus
\[  D=\{r e^{i \theta} : e^{-n} < r < 1, 0 < \theta < 2 \pi \}, \]
with
\[  V_1 = \{|z| = e^{-n}\} , \;\; V_2 =
    \{|z| = 1 \} , \]
then $D$ is conformally equivalent to a rectangle of sides $n$ and $2 \pi$
and hence
\[  \Delta(D,V_1,V_2) = \frac{n}{2 \pi} . \]
Suppose $V_2$ is a closed set such that every  $\gamma:[0,1]
\rightarrow \C$ with $\gamma(0) \in V_1, \gamma(1) \in V_3,
\gamma(0,1) \subset D$ satisfies
\[  \gamma(0,1) \cap V_2 \neq \emptyset . \]
Then \cite[Proposition 9.2]{Pom}
\[    \Delta(D,V_1,V_3) \geq \Delta(D,V_1,V_2) + \Delta(D,V_2,V_3) . \]
In other words, the extremal length satisfies a reversed triangle
inequality.  In particular if we let $\disk$ denote the unit disk
and $\partial_n = \{|z| = e^{-n}\}$, then for any $D \subset \disk$,
\begin{equation}  \label{extremal.2}
   \Delta(D,\partial_n,\partial_0) \geq \sum_{j=1}^n 
     \Delta(D,\partial_j,\partial_{j-1}) . 
\end{equation}

Let $E$ be an interval in the unit circle
\[  E = \{e^{i \theta} : \theta_1 \leq \theta \leq \theta_2 \} , \]
with length $l(E) = \theta_2 - \theta_1$.
Let $\partial_r$ denote the circle of radius $e^{-r}$ as above,
and let $\Gamma_1 = \Gamma_1(r,E) = \Gamma(\disk,\partial_r,E)$. Let
$B_t$ denote a standard two dimensional Brownian motion, considered as
a complex valued Brownian motion, with stopping times
\[   T_n = \inf\{t: |B_t| = e^n\} . \]
 By
Pfluger's Theorem \cite[Theorem 9.17]{Pom}, for $1 \leq r \leq 5$,
\[  l(E) = \frac{1}{2 \pi} \Prob^0\{B(T_0) \in E \} 
    \asymp \exp\{- \pi \Delta_1 \} , \]
where $\Delta_1 = \Delta_1(\disk,\partial_r,E)$.  Let $\Gamma_2
= \Gamma_2(r)$ be the collection of curves $\{\gamma_\theta\}$,
\[  \gamma_\theta(t) = e^{- i \theta} [(1-t) e^{-r-1} + t e^{-r}], \;\; 
   0 \leq t \leq 1 . \]
Let $\Gamma_3 = \Gamma_3(r,E)$ denote the set of curves obtained
by combining any $\gamma_2 \in \Gamma_2$ with any $\gamma_1 \in \Gamma_1$.
The way to combine is to start with $\gamma_2;$ then take a curve
in $\partial_r$ which goes from $\gamma_2(1)$ to $\gamma_1(0)$; and
then follow $\gamma_1$.  Of course, we must do a simple reparameterization
to get a curve $\gamma:[0,1] \rightarrow \C$.  It is not difficult
to show that there is a $c$, independent of $r$ and $E$, such that
\[  \Delta(\Gamma_3) \leq \Delta_1 + c . \]
(In proving this one notes that one can restrict the set of $\rho$ to those
that are zero on $\partial_r$.)
Similarly, let $\Gamma_4$ be the collection of curves $\{\gamma_\theta\}$,
\[  \gamma_\theta(t) = e^{- i (\theta + 2\pi t)} [(1-t) e^{-r-1} + t e^{-r}] ,
  \; \;  0 \leq t \leq 1 , \]
and let $\Gamma_5$ be the collection of curves obtained by attaching 
curves from $\Gamma_4$ with curves from $\Gamma_1$.  Again we can show
that 
\[  \Delta(\Gamma_5) \leq \Delta_1 + c . \]
By comparison with $\Gamma_4$ and $\Gamma_5$, we can see that if
$h:[0,1] \rightarrow \C$ is any continuous 
curve with $h(0) \in \partial_{r+1}, h(1) \in
\partial_{r}$, then
\begin{equation}  \label{extremal.3}
  \Delta(\disk,h[0,1],E) \leq \Delta_1 + c . 
\end{equation}

Now suppose $r \geq 2$ and
$h:[0,1] \rightarrow \C$ is a continuous curve without double points
 with
$h(0) \in \partial_r$, $h(1) \in \partial_0$, $h[0,1) \subset \disk$.
Let $D$ be the connected component of $\disk \setminus h[0,1]$
whose boundary includes $\partial_0$.  Let $x \in \partial_r \cap D$
(if $\partial_r \cap D = \emptyset$, the conclusion we will obtain
below is trivial).  Let $d$ be the distance from $x$ to $\partial D$,
and note that $d \leq 2|x|$; in particular, the closed disk of radius
$d$ about $x$ is contained in the disk of radius $e^{-r+2}$ about
the origin.  Let $f$ be a conformal transformation of $D$ to the
unit disk with $f(x) = 0$.  Let $\delta$ be the smallest number such
that $f(\bar{{\cal B}}(x,\delta)) \cap \partial_2 \neq \emptyset$, where
$\bar{{\cal B}}(x,\delta)$ denotes the closed ball of radius $\delta$
about $x$.   Note
that $\delta < d$. By conformal invariance of Brownian motion (or
equivalently, by conformal invariance of harmonic measure),
\[  \Prob^x\{B[0,T_0] \cap h[0,1] = \emptyset\}
    = \Prob^0\{B(T_0) \in f(\partial_0)\} . \]
By Pfluger's Theorem, the right hand side is comparable to
\[  \exp\{- \pi \Delta(\disk,\partial_2,f(\partial_0))\} . \]
By (\ref{extremal.3}), this is comparable to
\[  \exp\{ - \pi \Delta(\disk,f(\bar{{\cal B}}(x,\delta)),f(\partial_0)) \}, \]
which by conformal invariance of extremal distance equals
\[  \exp\{ - \pi \Delta(D,\partial {\cal B}(x,\delta),\partial_0)\} . \]
Since ${\cal B}(x,\delta) \subset {\cal B}(0,e^{r-2})$, this last term
is smaller than
\[ \exp\{-\pi \Delta(D,\partial_{r-2},\partial_0)\} . \]
Combining all of this we get that there is a constant $c$ such that
if $h$ is as above and $x \in \partial_r$,
\begin{equation}  \label{extremal.4}
\Prob^x\{B[0,T_0] \cap h[0,1] = \emptyset \} 
   \leq c \exp\{-\pi \Delta(D,\partial_{r-2},\partial_0)\} . 
\end{equation}

Let $h,D$ be as in the previous paragraph,
and let $A$ be the connected component of $D \cap \{|z| > e^{-r} \}$
that contains the unit circle in its boundary.  Assume that
$\partial A \cap \partial_r \neq \emptyset$.  Note that if $\gamma$ is
a curve with $\gamma(0) \in \partial_r, \gamma(1) \in \partial_0, 
\gamma(0,1) \subset D$, then there exists an $s \geq 0$ such that
$\gamma(s) \in \partial_r, \gamma(s,1) \subset A$.  Hence,
\[  \Delta(D,\partial_{r-2},\partial_0) = \Delta(A,\partial_{r-2},\partial_0)
  . \]
It will be convenient if we conformally map this region by a logarithm.
In this case $\partial_r$ is sent to $U_r$ where
\[   U_j = \{\Re(z) = -j \} ; \]
We take an $h:[0,1] \rightarrow \C$ satisfying
\[  h(0) \in U_r, \;\; h(1) \in U_0,\;\; h(0,1) \subset \{\Re(z) < 0 \} , \]
and for each integer $k$ we have the curve
\[   h_k(t) = h(t) + 2\pi i , \;\;\; 0 \leq t \leq 1 . \]
Let $G$ be the region bounded by $U_r,U_{0}, h[\sigma_r,1]$
and $h_1[\sigma_r,1]$. 
The region $A$ is the connected component of $G$ whose boundary intersects
$U_0$.  Note that any curve connecting $U_r$ to $U_0$ in $G$ has a
subpath connecting $U_r$ to $U_0$ in $A$.  Since $A \subset G$,
\[   {\rm area}(A) \leq {\rm area}(G) \leq  2 \pi r . \]
If $j$ is a positive integer, let $A_j$ denote the connected component
of $A \cap \{\Re(z) < -j\}$ whose boundary includes a portion of
$U_r$, and let $V_j = U_j \cap \partial A_j$.  Note that $V_j$ is contained
in an interval of length at most $2 \pi$, and that every continuous curve
from $U_r$ to $U_0$ staying in $A$ hits $U_j$ first in $V_j$.
Fix some $0 < \epsilon < 1/10$.  Suppose we can find $k$ integers
$1 \leq j_1 < \cdots < j_k \leq r $ such that the following holds. 
For each $j \in \{j_1,\ldots,j_k\}$,  there
exists a point $z_j$ with
\[         \Re(z_j) = -j + \frac{1}{2}  , \]
and such that for each $w_1,w_2 \in {\cal B}(z_j,\epsilon)$,
\begin{equation}  \label{extremal.5}
 \dist(w_1,V_j + 2 \pi l) + \dist(w_2,V_{j-1} + 2 \pi l) \geq 1 , 
\end{equation}
for every integer $l$. 
In other words the portion of every path in $A$
 from $V_j$ to  $V_{j-1}$ outside of the $2 \pi i$ integer translates
of ${\cal B}(z_j,\epsilon)$ must have
length at least one.  Note that
\[  {\rm area}[G \setminus \bigcup_{j=1}^k \bigcup_{l = -\infty}^\infty
     {\cal B}(z_j + 2 \pi l,\epsilon)] \leq 2 \pi r - k \pi \epsilon^2. \]
If we let $\rho$ be the function that is $1$ on
\[    A \setminus[ \bigcup_{j=1}^k \bigcup_{l = -\infty}^\infty
     {\cal B}(z_j + 2 \pi l,\epsilon)], \]
and zero elsewhere, we have
\[    a(\rho) \leq 2 \pi r - k\pi \epsilon^2, \]
and
\[  L(\rho) \geq r . \]
Hence
\[  \Delta(A,U_r,U_0) \geq \frac{r^2}{2 \pi r - k\pi \epsilon^2}
    \geq \frac{r}{\pi}[\frac{1}{2} + \frac{\epsilon^2 k}{4r} ]. \]

Now fix a $j$ and suppose $s<t$ with $h(s) \in U_j, h(t) \in
U_{j-1}, h(s,t) \subset \{\Re(z) < j-1\}$. Suppose also that
\[     \Im(h(t)) \geq \Im(h(s)) + \delta , \]
for some $\delta > 0$.  Note that
\[   V_j \subset \{-j + iy: h(s) \leq y \leq h(s) + 2 \pi \}, \]
\[   V_{j-1} \subset \{-j+1 + i y : h(t) \leq y \leq h(t) + 2 \pi\} . \]
Let
\[     z_j = -j+ \frac{1}{2} + [h(s) + \frac{\delta}{2}]i  . \]
 Note
 that if we draw a sufficiently small ball around this point
(we leave the high school geometry estimate to the reader), that
it satisfies (\ref{extremal.5}) for some $\epsilon = \epsilon_\delta$.
Similarly if $\Im(h(t)) \leq \Im(h(s)) - \delta$, a similar fact
holds using 
\[  z_j = -j+ \frac{1}{2} + [h(s) - \frac{\delta}{2}]i  . \]
When we transform this argument back to the unit disk and use
(\ref{extremal.4}), we get
the following lemma.  (The $\epsilon$ in this lemma corresponds to
the $\epsilon^2/4$ above.)

\begin{lemma}  \label{extlemma.1}
For every $\delta > 0 $, there exists an $\epsilon > 0$ and
a $c < \infty$ such that the following is true.  Let $r$ be a positive
integer and $h:[0,1] \rightarrow \C$ a continuous function without double
points with
$h(0) \in \partial_r, h(1) \in \partial_0,h(0,1) \subset \disk$.
For each integer $j=1,\ldots,r-2$, let 
\[ \nu_j = \inf\{t: h(t) \in \partial_j \} , \]
\[  Y_j = |\arg(h(\nu_j) ) - \arg(h(\nu_{j-1}))| . \]
Let
\[  W = W_r(h,\delta) = \sum_{j=1}^r I\{Y_j \geq \delta\} . \]
Then if $x \in \partial_r$,
\[ \Prob^x\{B[0,T_0] \cap h[0,1] = \emptyset \} 
   \leq c \exp\{-r(\frac{1}{2}) +  \epsilon W)\} . \]
\end{lemma}

\section{Strong Approximation} \label{strongsec}

We will need to use a strong approximation of a simple random walk
and a Brownian motion in two dimensions.  The approximation we will
use is derived from the standard Skorohod embedding of a one dimensional
simple random walk in a one dimensional Brownian motion.  Because the
result we need is slightly different than other results (in particular,
the independence of the event $U_n$ below from a particular $\sigma$-algebra),
we will sketch the derivation.  Let $(\Omega,{\cal F},\Prob)$ be a probability
space on which are defined
 a two-dimensional Brownian motion, $B_t = B(t) = (B_t^1,B_t^2)$,
and a one-dimensional simple random walk, $W_k$, that is independent
of the Brownian motion.  As before, let
\[  T_n = \inf\{t: |B_t| = e^n \}, \]
and let
\[   R_k = \frac{1}{2}(W_k + k) . \]

Let $S_k^i$ be the simple random walk derived from $B_t^i$ from the
Skorohod embedding.  This is obtained by setting $\eta_0^i = 0,$
\[  \eta_{k+1}^i = \inf\{t \geq \eta_k^i: |B^i(t) - B^i(\eta_{k}^i)|
   = 1 \}, \]
and
\[  S^i_k = B^i(\eta_k^i) . \]
Since $\E[\eta_1^i] = 1$,
one expects $\eta_k^i - k$ to be of order $k^{1/2}$ and hence
$B(\eta_k^i) - B(k)$ to be of order $k^{1/4}$.  Standard techniques make
this precise; in particular, since $\eta_k^i$ has an exponential moment, we
can derive exponential estimates on the
probabilities.  In fact, one can show that there exists $c,\beta$ such that
\[  \Prob \{ \; \sup_{0 \leq s,t \leq e^{9n/8}, |s-t| \leq e^{5n/8} }
    |B_t^i - B_s^i| \geq e^{3n/8} \; \} \leq c \exp\{-e^{\beta n} \} . \]
(In this section we let $c,\beta$ be positive constants whose value
may change from line to line.)  
Also,
\[  \Prob \{ \; \sup_{0 \leq k \leq e^{9n/8}} |\eta^i_k - k |
      \geq e^{5n/8} \; \} \leq c \exp\{-e^{\beta n} \} . \]
In particular,
\[  \Prob\{\sup_{0 \leq t \leq e^{9n/8}} |B^i_t - S^i( \lf t \rf)|
     \geq e^{3n/8} \} \leq c  \exp\{-e^{\beta n} \} . \]

If we let
\[   S_k = (S^1(R_k),S^2(k - R_k)) ,\]
it is easy to see that $S_k$ is a simple random walk in $\Z^2$.
Again, standard estimates given
\[       \Prob\{ \; \sup_{0 \leq k \leq e^{9n/8}}
      |R_k - \frac{k}{2}| \geq e^{5n/8} \; \} \leq c 
    \exp\{-e^{\beta n} \} , \]
and hence
\[  \Prob\{\sup_{0 \leq t \leq e^{9n/8}} |B_t - S(\lf 2t \rf)|
     \geq c  e^{3n/8} \} \leq c  \exp\{-e^{\beta n} \} . \]
Also note that
\[  \Prob\{T_n \geq e^{9n/4} \} \leq c  \exp\{-e^{\beta n} \} . \]

We now let $U_n$ be the event
\[  U_n = \{T_n \leq e^{9n/4}; \sup_{0 \leq t \leq T_n} |B_t - S(\lf 2t \rf )|
     \leq c  e^{3n/4} \}. \]
This event is measurable with respect to the $\sigma$-algebra generated
by
\[ \{B_t: t \leq T_n\} \cup \{R_k,k=0,1,2,\ldots\}. \]
In particular it is independent of the $\sigma$-algebra generated
by
\[   \{B(t+T_n) - B(T_n): t \geq 0 \} . \]
We have sketched the proof of the following. The Brownian motion $B_t$ in
the lemma has variance parameter $1/2$, i.e., $B_t = \tilde{B}_{t/2}$,
where $\tilde{B}$ is a standard Brownian motion. 

\begin{lemma}  \label{lemma.strong}
There exists a $c,\beta$ such that a planar Brownian motion with
variance parameter $1/2$
and a two-dimensional simple random walk $S_k$ can
be defined on the same probability space satisfying the following.
 For each $n$, 
there exists an
event $E_n$ that is independent of
\[  \{B(t + T_n) - B(T_n): t \geq 0 \}, \]
with
\[  \Prob(E_n) \geq 1 - c \exp\{-e^{\beta n}\} , \]
and such that on the event $E_n$,
\[   |B(t) - S(\lf t \rf )| \leq e^{7n/8} , \;\;\; t \leq T_{n-(1/2)} . \]
\end{lemma}

\vspace{2ex}

Let 
\[  \tau_n = \inf\{k: |S_k| = e^n \} . \]
One can see that if $V$ is any event that is measurable with respect
to
\[  \{S_k : 0 \leq k \leq \tau_{n-(1/2)} \}, \]
then $V \cap E_n$ is also independent of the $\sigma$-algebra generated
by
\[   \{B(t+T_n) - B(T_n): t \geq 0 \} . \]

\section{Bound for Random Walks}  \label{walksec}

Let $ \delta_0,\epsilon$ be such that Lemma
\ref{extlemma.1} holds (for some constant $c$), and choose
$\delta > \delta_0$.  We allow constants $c,c_1,c_2,\ldots$ in this section
to depend on $\delta_0,\delta,\epsilon$. The goal of this section is
to show that a corresponding result holds for simple random walk
for $\delta,\epsilon$.  We first use Brownian scaling to
make a slight restatement of Lemma \ref{extlemma.1}.  If $h:[0,\infty)
\rightarrow \C$ is a continuous curve without double points
 with $h(0) = 0$, $|h(t)|
\rightarrow \infty, t \rightarrow \infty$, let
\[  \sigma_n = \inf\{t: |h(t)| = e^n \} , \]
\[ Y_n = Y_n(h) = |\arg(h(\sigma_n)) - \arg(h(\sigma_{n+1}))| , \]
\[W_{m,n} = W_{m,n}(h,\delta_0) = \sum_{k=0}^{n-1}
     I\{Y_{m+k} \geq \delta_0\} . \]
  Then Lemma
\ref{extlemma.1} immediately implies the following.

\begin{lemma}  \label{walklem.1}
There exists a constant $c$ such that for all nonnegative integers
$m,n$ and all $h$, 
\[  \Prob^0\{B[T_m,T_{m+n}] \cap h[\sigma_m,\sigma_{m+n}]
    = \emptyset\} \leq c \exp\{-n\frac{1}{2} - \epsilon W_{m,n} \} . \]
\end{lemma}

Let $\omega:\{0,1,2,\ldots\} \rightarrow \Z^2$ be a self-avoiding
random walk path with $\omega(0) = 0, |\omega(t)| \rightarrow \infty,\:
t \rightarrow \infty$.  Associated with $\omega$ is the continuous path
$h_\omega: [0,\infty) \rightarrow \C$ obtained by linear interpolation,
i.e., 
\[   h_\omega(t) = \omega(\lfloor t \rfloor) + (t - \lfloor t \rfloor)
       [\omega(\lfloor t \rfloor + 1) - \omega(\lfloor t \rfloor)]\]
(where we consider $\Z^2$ as embedded in $\C$).  Let
\[  \tilde{\sigma}_n = \inf\{t: |\omega(t)| \geq e^{n} \} , \]
\[    \tilde{Y}_n = \tilde{Y}_n(\omega) = 
    |\arg(\omega(\tilde{\sigma}_n)) -
   \arg(\omega(\tilde{\sigma}_{n+1})) | , \]
\[  \tilde{W}_{m,n} = \tilde{W}_{m,n}(\omega,\delta) =
\sum_{k=0}^{n-1}
     I\{\tilde{Y}_{m+n} \geq \delta\} . \]
Note that if  $\sigma_n$ are the times defined as above for $h = h_\omega$,
\[      \sigma_n \leq \tilde{\sigma}_{n} \leq \sigma_n + 1 . \]
In particular, for all $n$ sufficiently large, if $\tilde{Y}_n \geq
\delta$, then $Y_n \geq \delta_0$. Therefore,
\begin{equation}  \label{walk.3}
   W_{m,n}(h_\omega) \geq \tilde{W}_{m,n}(\omega) - c_1 . 
\end{equation}
 Our goal is to prove the following.

\begin{lemma}  \label{walklem.2}
 There exists a constant $c_2$
such that for all nonnegative integers
$m,n$, and all $\omega$ as above,
\[  \Prob^0\{S[\tau_m,\tau_{m+n}] \cap \omega[\tilde{\sigma}_m,
  \tilde{\sigma}_{m+n}]
    = \emptyset\} \leq c_2 \exp\{-n\frac{1}{2} - \epsilon
   \tilde{W}_{m,n} \} . \]
\end{lemma}

Let
\[  b_n = b_n(\delta,\epsilon) =
  \sup \exp\{n \frac{1}{2} + \epsilon \tilde{W}_{m,n} \}  \Prob^0\{S[\tau_m,\tau_{m+n}] \cap \omega[\tilde{\sigma}_m,
  \tilde{\sigma}_{m+n}]
    = \emptyset\} , \]
where the supremum is over all $\omega$ and all positive integers $m$.
Lemma \ref{walklem.2} is equivalent to saying that the sequence
$\{b_n\}$ is bounded.  It is obvious that $b_0 = 1$ and since
$\tilde{W}_{m,n} \leq n$, 
\[  b_{m+n} \leq e^{n \alpha} b_m , \]
where $\alpha = (1/2) + \epsilon$.
 If we show that there
is a $\beta> 0$, and a $c_3 < \infty$ such that for all $n$,
\begin{equation}  \label{walk.1}
  b_n \leq c_3 \sum_{j=0}^{n-1} e^{- j \beta } b_j , 
\end{equation} 
then it follows that the $b_n$ are bounded (see, for example,
\cite[Lemma 4.5]{walkcut}).  So to prove Lemma \ref{walklem.2}
it suffices to prove (\ref{walk.1}) for $n \geq 3$.

Let $B$ be a Brownian motion starting at the origin in $\C$ and
let $S$ be the corresponding simple random walk derived from the
strong approximation as in Lemma \ref{lemma.strong}.  Let $\omega$
be a self-avoiding random walk path
and $h = h_\omega$ as above.  For 
$1 \leq k \leq n-1$, let $U_k = U_k(m,n)$ be the event
\[  U_k = \{B[T_{m+k-1},T_{m+k }] \cap h[\sigma_m,\sigma_{m+n}]
     \neq \emptyset; \hspace{1in} \]
\[  \hspace{1.5in} B[T_{m+k},T_{m+n-1}] \cap h[\sigma_m,\sigma_{m+n}]
   = \emptyset \}. \]
Let $U_0 = U_0(m,n)$ be the event
\[ U_0 = \{B[T_m,T_{m+n-1}] \cap  h[\sigma_m,\sigma_{m+n}]
  = \emptyset\} , \]
and let
 $V = V_{m,n}$ be the event
\[  \{S[\tau_{m},\tau_{m + n}] \cap \omega[\tilde{\sigma}_m,
  \tilde{\sigma}_{m+n}] = \emptyset \}. \]
Note that
\[  V = \bigcup_{k=0}^{n-1} (V \cap U_k) . \]
To prove (\ref{walk.1}) it suffices to prove that for some 
positive constants
$c,\beta$,
\begin{equation}  \label{walk.2}
    \Prob[V \cap U_k] \leq c b_k e^{- k \beta } \exp\{-\frac{1}{2}n
   - \epsilon \tilde{W}_{m,n}\} . 
\end{equation}
For $k=0,1$,  this estimate follows immediately from Lemma \ref{walklem.1}
and (\ref{walk.3}) and hence we will assume $k \geq 2$.  Note that
\[      \tilde{W}_{m,m+k-2} + W_{m+k+3,m+n} 
    \geq \tilde{W}_{m,n} - c_3 . \]

Fix $k$ and let $E = E_{m+k+2}$ be the event as in Lemma \ref{lemma.strong}
and write $\beta_1$ for the $\beta$ in that lemma.  Note that
\begin{align*}
\Prob[V \cap U_k \cap E^c] &\leq \Prob[E^c;
         B[T_{m+k+2},T_{m+n-1}] \cap h[\sigma_{m+k+2},\sigma_{m+n}]
    = \emptyset \} \\
    &\leq c  \Prob(E^c) \Prob \{B[T_{m+k+3},T_{m+n-1}] \cap h[\sigma_{m+k+2},\sigma_{m+n}]
    = \emptyset \} \\
  & \leq c \exp\{-e^{\beta_1 k}\}  \exp\{\frac{1}{2}(n-k) - \epsilon
      W_{m+k+3,m+n-1} \} \\
  & \leq c \exp\{-\frac{1}{2}n - \epsilon W_{m,m+n-1} \} .
\end{align*}
The second inequality uses the conditional independence and the
Harnack inequality for Brownian motion (harmonic functions), and
the last inequality uses the trivial inequality ${W}_{m,m+k}
    \leq k $.  Note that $W_{m,n+m-1} \geq \tilde{W}_{m,n} - c$.

Let $F = F_{m,k}$ be the event
\[ F = \{S[\tau_m,\tau_{m+k+1}] \cap \omega[\hat{\sigma}_m,
  \hat{\sigma}_{m+k+1}] = \emptyset; \dist(S[\tau_{m+k-2},\tau_{m+k}],
   \omega[\hat{\sigma}_{m},\hat{\sigma}_{m+k}])
    \leq e^{15(m+k)/16} \} , \]
and let
\[ \rho = \rho_{m,k,\omega} = \inf\{t \geq \tau_{m+k-2}:
       \dist(S(t),\omega[\tau_{m+k-2},\tau_{m+k}])
    \leq e^{15(m+k)/16} \} . \]
By the discrete Beurling projection theorem (see \cite[Lemma 2.3]{LMak}
for the more general version used here) and the strong Markov property
for random walk,
\[  \Prob\{S[\rho,\tau_{m+k+1}] \cap \omega[\hat{\sigma}_m,
  \hat{\sigma}_{m+k+1}] = \emptyset \mid \rho \leq
  \tau_{m+k}    \} \leq c e^{-k/32} . \]
Hence,
\[ \Prob(F) \leq c b_k \exp\{-k \frac{1}{2} - \epsilon \tilde{W}_{m,m+k-2}
  - k \frac{1}{32}\} . \]
However, conditioned on the event $E$, the event $F$ is independent of
\[ \{B(t + T_{m+k+2}) - B(T_{m+k+2}): t \geq 0 \} . \]
Also
\[  V \cap U_k \cap E \; \subset \; F \cap \{B[T_{m+k+2},T_{m+n-1}]
   \cap h[\sigma_{m+k+2},\sigma_{m+n-1}] = \emptyset \}. \]
Hence from the strong Markov property for Brownian motion and the Harnack
inequality,
\begin{align*}
 \Prob(V \cap U_k \cap E) &\leq c b_k e^{-k/32} \exp\{-n \frac{1}{2}
     - \epsilon (\tilde{W}_{m,m+k-2} + W_{m+k+2,m+n-1})\} \\
    &\leq
    c b_k e^{-k/32} \exp\{-n \frac{1}{2} - \epsilon \tilde{W}_{m+n} \},
\end{align*}
which gives (\ref{walk.2}) and hence proves Lemma \ref{walklem.1}.

\end{document}